\documentclass[preprint,12pt]{elsarticle}

\usepackage{amssymb}
\usepackage{amsmath}
\journal{Discrete Applied Mathematics}

\usepackage{amsthm, enumitem, bm, xcolor}
\theoremstyle{plain}
\newtheorem{theorem}{Theorem}
\newtheorem{lemma}[theorem]{Lemma}
\newtheorem{corollary}[theorem]{Corollary}
\theoremstyle{definition}
\newtheorem{definition}[theorem]{Definition}

\usepackage{subcaption}
\usepackage{array}
\usepackage{longtable}

\usepackage{tikz}
\usetikzlibrary{arrows.meta}
\usetikzlibrary{automata}

\usepackage[hyphens]{url}
\usepackage[hidelinks]{hyperref}

\usepackage{afterpage}

\begin{document}

\begin{frontmatter}

%% Title, authors and addresses

%% use the tnoteref command within \title for footnotes;
%% use the tnotetext command for theassociated footnote;
%% use the fnref command within \author or \affiliation for footnotes;
%% use the fntext command for theassociated footnote;
%% use the corref command within \author for corresponding author footnotes;
%% use the cortext command for theassociated footnote;
%% use the ead command for the email address,
%% and the form \ead[url] for the home page:
%% \title{Title\tnoteref{label1}}
%% \tnotetext[label1]{}
%% \author{Name\corref{cor1}\fnref{label2}}
%% \ead{email address}
%% \ead[url]{home page}
%% \fntext[label2]{}
%% \cortext[cor1]{}
%% \affiliation{organization={},
%%             addressline={},
%%             city={},
%%             postcode={},
%%             state={},
%%             country={}}
%% \fntext[label3]{}

\title{Classification of borderenergetic chemical graphs \\ and borderenergetic graphs of order~$12$}

%% use optional labels to link authors explicitly to addresses:
%% \author[label1,label2]{}
%% \affiliation[label1]{organization={},
%%             addressline={},
%%             city={},
%%             postcode={},
%%             state={},
%%             country={}}
%%
%% \affiliation[label2]{organization={},
%%             addressline={},
%%             city={},
%%             postcode={},
%%             state={},
%%             country={}}

\author[1a,1b]{Péter Csikvári}
\author[2,3,4]{Ivan Damnjanović\corref{cor1}}
\author[5]{Marko Milošević}
\author[5]{Ivan Stanković}
\author[6]{Dragan Stevanović}

\cortext[cor1]{Corresponding author.}

%% Author affiliation
\affiliation[1a]{
    organization={HUN-REN Alfréd Rényi Institute of Mathematics},
    addressline={Reáltanoda utca 13--15},
    city={Budapest},
    postcode={1053},
    country={Hungary}
}

\affiliation[1b]{
    organization={Eötvös Loránd University},
    addressline={Pázmány Péter sétány 1/C},
    city={Budapest},
    postcode={1117},
    country={Hungary}
}

\affiliation[2]{
    organization={Faculty of Mathematics, Natural Sciences and Information Technologies, \\ University of Primorska},
    addressline={Glagoljaška 8},
    city={Koper},
    postcode={6000},
    country={Slovenia}
}

\affiliation[3]{
    organization={Faculty of Electronic Engineering, University of Niš},
    addressline={\\ Aleksandra Medvedeva 4},
    city={Niš},
    postcode={18104},
    country={Serbia}
}

\affiliation[4]{
    organization={Diffine LLC},
    addressline={3681 Villa Terrace},
    city={San Diego},
    postcode={92104},
    state={California},
    country={USA}
}

\affiliation[5]{
    organization={Faculty of Sciences and Mathematics, University of Niš},
    addressline={\\ Višegradska 33},
    city={Niš},
    postcode={18106},
    country={Serbia}
}

\affiliation[6]{
    organization={College of Integrative Studies, Abdullah Al-Salem University},
    addressline={\\ Firdous Street, Block~3},
    city={Khaldiya},
    postcode={72303},
    country={Kuwait}
}

%\affiliation[7]{
%    organization={Mathematical Institute, Serbian Academy of Sciences and Arts},
%    addressline={Kneza Mihaila 36},
%    city={Belgrade},
%    postcode={11001},
%    country={Serbia}
%}

%% Abstract
\begin{abstract}
The energy $\mathcal{E}(G)$ of a simple graph~$G$ is the sum of absolute values of the eigenvalues of its adjacency matrix. 
A borderenergetic graph of order $n \in \mathbb{N}$ is any noncomplete graph~$G$ 
such that $\mathcal{E}(G) = \mathcal{E}(K_n) = 2n - 2$.
Here we combine two-phase computer-assisted search with theoretical arguments to show that 
there are only three borderenergetic chemical graphs,
thus completing the earlier findings of Li, Wei and Zhu [MATCH Commun.\ Math.\ Comput.\ Chem.\ 77 (2017), 25--36].
We perform two-phase computer-assisted search to also find all $566$ borderenergetic graphs of order~$12$,
thereby correcting and extending the results from a previous search performed by Furtula and Gutman 
[Iranian J.\ Math.\ Chem.\ 8(4) (2017), 339--344]. 
\end{abstract}


\end{frontmatter}

%% Add \usepackage{lineno} before \begin{document} and uncomment 
%% following line to enable line numbers
%% \linenumbers

%% main text
%%

%% Use \section commands to start a section
\section{Introduction}

Introduced by Gutman~\cite{Gutman1978}, 
the energy $\mathcal{E}(G)$ of a simple graph $G$ is the sum of absolute values of the eigenvalues of its adjacency matrix. 
The graph energy has a wide application in chemical graph theory and 
has been extensively studied during the last two decades (see, e.g., \cite{GutRam2020, LiShiGut2012}). 
Here, we will deal with the \emph{borderenergetic graphs}, which are defined as follows.

\begin{definition}
A \emph{borderenergetic graph} of order $n \in \mathbb{N}$ is a noncomplete graph~$G$ 
such that $\mathcal{E}(G) = \mathcal{E}(K_n) = 2n - 2$.
\end{definition}

Many graph theorists have found interest in investigating the borderenergetic graphs 
(see, e.g., the survey paper \cite{GhorDengHakLi2020}, and the references therein). 
In a series of papers~\cite{GongLiXuGutFur2015, LiWeiGong2015, ShaoDeng2016}
it was proved that there is no borderenergetic graph of order below seven, 
while for each $n \ge 7$ there necessarily exists a borderenergetic graph of order~$n$. 
The numbers of borderenergetic graphs of order~$n$ were determined for $7 \le n \le 11$
in~\cite{GongLiXuGutFur2015,LiWeiGong2015,ShaoDeng2016},
as shown in Table~\ref{table:borderenergetic_count}. 
In a subsequent paper, 
Furtula and Gutman \cite[Theorem 2]{FurGut2017} performed a computer-assisted search 
that incorrectly showed the existence of $572$ connected borderenergetic graphs of order~$12$.

\begin{table}[h!]
\centering
\begin{tabular}{|l|*{7}{c|}}\hline
Order & 1 -- 6 & 7 & 8 & 9 & 10 & 11 & 12 \\\hline
Borderenergetic graphs & $0$ & $1$ & $6$ & $17$ & $49$ & $158$ & $566$ \\\hline
\end{tabular}
\caption{The number of borderenergetic graphs of orders up to~$12$.}
\label{table:borderenergetic_count}
\end{table}

Some researchers have analyzed borderenergetic graphs that belong to specific classes. 
A \emph{chemical graph} is a graph whose maximum degree is at most four.
The following two theorems of Li, Wei and Zhu tell us that 
a borderenergetic chemical graph has at most $21$ vertices.

\begin{theorem}[\hspace{1sp}{\cite[Theorem 3.1]{LiWeiZhu2017}}]\label{li_th_1}
There is no borderenergetic graph with maximum degree two or three.
\end{theorem}

\begin{theorem}[\hspace{1sp}{\cite[Theorem 3.2]{LiWeiZhu2017}}]\label{li_th_2}
    Let $G$ be a borderenergetic graph of order~$n$ with maximum degree four. 
    Then $G$ must have the following properties:
    \begin{enumerate}[label=\textbf{(\roman*)}]
        \item $|E(G)| = 2n$ or $|E(G)| = 2n - 1$;
        \item $n \le 21$;
        \item $G$ is nonbipartite;
        \item the nullity, i.e., the multiplicity of eigenvalue zero, of $G$ is zero.
    \end{enumerate}
\end{theorem}

Here we investigate the borderenergetic chemical graphs and the borderenergetic graphs of order~$12$. 
We apply two-phase computer-assisted search
together with a mathematical optimization theory based proof 
to provide the full classification of borderenergetic chemical graphs,
thus solving the problem initiated in~\cite{LiWeiZhu2017}.

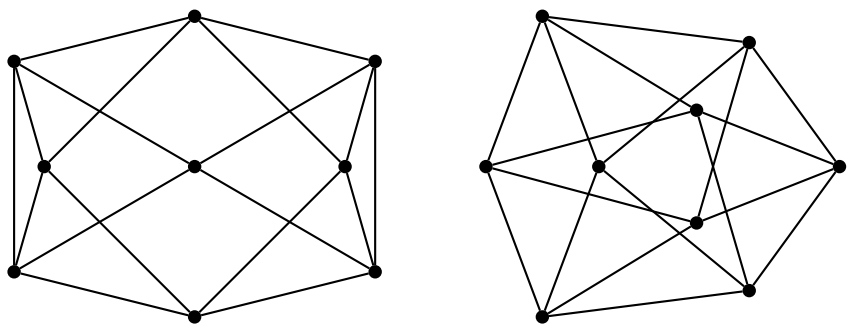
\begin{figure}[htbp]
    \centering
    \hspace{1.0cm}
    \begin{subfigure}[b]{0.38\textwidth}
        \centering
        \begin{tikzpicture}[scale=1.0]
            \node[circle, draw=black, fill=black, inner sep=0pt, minimum size=0.15cm, thick] (1) at (0, 0) {};
            \node[circle, draw=black, fill=black, inner sep=0pt, minimum size=0.15cm, thick] (2) at (0, 2) {};
            \node[circle, draw=black, fill=black, inner sep=0pt, minimum size=0.15cm, thick] (3) at (0, -2) {};
    
            \node[circle, draw=black, fill=black, inner sep=0pt, minimum size=0.15cm, thick] (4) at (2, 0) {};
            \node[circle, draw=black, fill=black, inner sep=0pt, minimum size=0.15cm, thick] (5) at (2.4, 1.4) {};
            \node[circle, draw=black, fill=black, inner sep=0pt, minimum size=0.15cm, thick] (6) at (2.4, -1.4) {};
    
            \node[circle, draw=black, fill=black, inner sep=0pt, minimum size=0.15cm, thick] (7) at (-2, 0) {};
            \node[circle, draw=black, fill=black, inner sep=0pt, minimum size=0.15cm, thick] (8) at (-2.4, 1.4) {};
            \node[circle, draw=black, fill=black, inner sep=0pt, minimum size=0.15cm, thick] (9) at (-2.4, -1.4) {};
    
            \draw[thick] (1) -- (5);
            \draw[thick] (1) -- (6);
            \draw[thick] (1) -- (8);
            \draw[thick] (1) -- (9);
    
            \draw[thick] (2) -- (4);
            \draw[thick] (2) -- (5);
            \draw[thick] (2) -- (7);
            \draw[thick] (2) -- (8);
    
            \draw[thick] (3) -- (4);
            \draw[thick] (3) -- (6);
            \draw[thick] (3) -- (7);
            \draw[thick] (3) -- (9);
    
            \draw[thick] (4) -- (5);
            \draw[thick] (5) -- (6);
            \draw[thick] (6) -- (4);
            \draw[thick] (7) -- (8);
            \draw[thick] (8) -- (9);
            \draw[thick] (9) -- (7);
        \end{tikzpicture}
    \end{subfigure}
    \hfill
    \begin{subfigure}[b]{0.38\textwidth}
        \centering
        \begin{tikzpicture}[scale=1.0]
            \node[circle, draw=black, fill=black, inner sep=0pt, minimum size=0.15cm, thick] (1) at (0, 0) {};
            \node[circle, draw=black, fill=black, inner sep=0pt, minimum size=0.15cm, thick] (2) at (-1.5, 0) {};        
            \node[circle, draw=black, fill=black, inner sep=0pt, minimum size=0.15cm, thick] (3) at (-0.75, 2) {};
            \node[circle, draw=black, fill=black, inner sep=0pt, minimum size=0.15cm, thick] (4) at (-0.75, -2) {};
    
            \node[circle, draw=black, fill=black, inner sep=0pt, minimum size=0.15cm, thick] (5) at (1.3, 0.75) {};
            \node[circle, draw=black, fill=black, inner sep=0pt, minimum size=0.15cm, thick] (6) at (1.3, -0.75) {};
            \node[circle, draw=black, fill=black, inner sep=0pt, minimum size=0.15cm, thick] (7) at (2, 1.65) {};
            \node[circle, draw=black, fill=black, inner sep=0pt, minimum size=0.15cm, thick] (8) at (2, -1.65) {};
    
            \node[circle, draw=black, fill=black, inner sep=0pt, minimum size=0.15cm, thick] (9) at (3.2, 0) {};
    
            \draw[thick] (1) -- (3);
            \draw[thick] (1) -- (4);
            \draw[thick] (2) -- (3);
            \draw[thick] (2) -- (4);
    
            \draw[thick] (5) -- (2);
            \draw[thick] (5) -- (3);
            \draw[thick] (5) -- (8);
            
            \draw[thick] (6) -- (2);
            \draw[thick] (6) -- (4);
            \draw[thick] (6) -- (7);
    
            \draw[thick] (7) -- (1);
            \draw[thick] (7) -- (3);
            \draw[thick] (8) -- (1);
            \draw[thick] (8) -- (4);
            
            \draw[thick] (9) -- (5);
            \draw[thick] (9) -- (6);
            \draw[thick] (9) -- (7);
            \draw[thick] (9) -- (8);
        \end{tikzpicture}
    \end{subfigure}
    \hspace{1.0cm}
    \caption{The two known borderenergetic chemical graphs on nine vertices~\cite{LiWeiZhu2017}.}
    \label{fig:chemical9}
\end{figure}

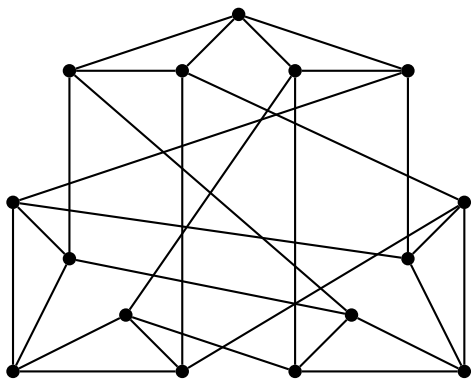
\begin{figure}[htbp]
    \centering
    \begin{tikzpicture}
        \node[circle, draw=black, fill=black, inner sep=0pt, minimum size=0.15cm, thick] (0) at (0, 4) {};
        \node[circle, draw=black, fill=black, inner sep=0pt, minimum size=0.15cm, thick] (1) at (-2.25, 3.25) {};
        \node[circle, draw=black, fill=black, inner sep=0pt, minimum size=0.15cm, thick] (2) at (-0.75, 3.25) {};
        \node[circle, draw=black, fill=black, inner sep=0pt, minimum size=0.15cm, thick] (3) at (0.75, 3.25) {};
        \node[circle, draw=black, fill=black, inner sep=0pt, minimum size=0.15cm, thick] (4) at (2.25, 3.25) {};

        \node[circle, draw=black, fill=black, inner sep=0pt, minimum size=0.15cm, thick] (5) at (-3, 1.5) {};
        \node[circle, draw=black, fill=black, inner sep=0pt, minimum size=0.15cm, thick] (6) at (-2.25, 0.75) {};
        \node[circle, draw=black, fill=black, inner sep=0pt, minimum size=0.15cm, thick] (7) at (3, 1.5) {};
        \node[circle, draw=black, fill=black, inner sep=0pt, minimum size=0.15cm, thick] (8) at (2.25, 0.75) {};

        \node[circle, draw=black, fill=black, inner sep=0pt, minimum size=0.15cm, thick] (9) at (-1.5, 0) {};
        \node[circle, draw=black, fill=black, inner sep=0pt, minimum size=0.15cm, thick] (10) at (-0.75, -0.75) {};
        \node[circle, draw=black, fill=black, inner sep=0pt, minimum size=0.15cm, thick] (11) at (1.5, 0) {};
        \node[circle, draw=black, fill=black, inner sep=0pt, minimum size=0.15cm, thick] (12) at (0.75, -0.75) {};

        \node[circle, draw=black, fill=black, inner sep=0pt, minimum size=0.15cm, thick] (13) at (-3, -0.75) {};
        \node[circle, draw=black, fill=black, inner sep=0pt, minimum size=0.15cm, thick] (14) at (3, -0.75) {};

        \draw[thick] (0) -- (1);
        \draw[thick] (0) -- (2);
        \draw[thick] (0) -- (3);
        \draw[thick] (0) -- (4);
        \draw[thick] (1) -- (2);
        \draw[thick] (3) -- (4);

        \draw[thick] (13) -- (5);
        \draw[thick] (13) -- (6);
        \draw[thick] (13) -- (9);
        \draw[thick] (13) -- (10);
        \draw[thick] (14) -- (7);
        \draw[thick] (14) -- (8);
        \draw[thick] (14) -- (11);
        \draw[thick] (14) -- (12);

        \draw[thick] (1) -- (6);
        \draw[thick] (1) -- (11);
        \draw[thick] (4) -- (5);
        \draw[thick] (4) -- (8);
        \draw[thick] (2) -- (7);
        \draw[thick] (2) -- (10);
        \draw[thick] (3) -- (9);
        \draw[thick] (3) -- (12);

        \draw[thick] (5) -- (6);
        \draw[thick] (7) -- (8);
        \draw[thick] (9) -- (10);
        \draw[thick] (11) -- (12);

        \draw[thick] (5) -- (8);
        \draw[thick] (7) -- (10);
        \draw[thick] (6) -- (11);
        \draw[thick] (9) -- (12);
    \end{tikzpicture}
    \caption{The newly found borderenergetic chemical graph on $15$ vertices.}
    \label{fig:chemical15}
\end{figure}

\begin{theorem}\label{main_chemical_th}
There exist exactly three borderenergetic chemical graphs, shown in Figs. \ref{fig:chemical9} and~\ref{fig:chemical15}.
\end{theorem}

We also perform a two-phase computer-assisted search that yields all borderenergetic graphs of order~$12$.
The obtained search results can be stated as follows.

\begin{theorem}\label{borderenergetic_12_th}
There exist precisely $566$ borderenergetic graphs of order~$12$, 
among which exactly seven are disconnected. 
The list of these graphs is given in~\cite{GitHubStuff}.
\end{theorem}

As a direct consequence, 
we conclude that the search results given in~\cite{FurGut2017} are not entirely correct 
due to the presence of numerical precision errors that occurred for $13$~connected graphs. 

Throughout the paper, we assume that all the graphs are undirected, finite and simple.
In Section~\ref{sc_computer_search} we elaborate the logic behind the two-phase computer-assisted search 
that was applied to find the borderenergetic graphs of order~$12$, 
as well as the borderenergetic chemical graphs on at most $19$~vertices. 
Subsequently, Section \ref{sc_nonexistence} demonstrates that 
there are no borderenergetic chemical graphs on $20$ or~$21$ vertices,
thereby completing the proof of Theorem~\ref{main_chemical_th}. 
Finally, in Section~\ref{sc_conclusion} we give a brief conclusion concerning all the obtained results.

\section{Two-phase computer-assisted search}\label{sc_computer_search}

We performed the computer-assisted search 
to find the borderenergetic graphs of order~$12$ and the borderenergetic chemical graphs on at most $19$~vertices
in two phases.
In the first phase we used the ``coarse sieve''
to select from the initial set of graphs
all graphs whose computed energy differs from $\mathcal{E}(K_n)=2n-2$ 
by a relatively large value ($10^{-6}$ to $10^{-8}$) 
when compared to the numerical precision ($10^{-12}$ to $10^{-15}$) 
offered by fast eigenvalue computation numerical routines.
This way we ensured that the selected candidates for borderenergetic graphs
necessarily include all true positives,
together with a reasonably small number of false positives.
In the second phase we used the ``fine sieve'',
by employing either symbolic computation in Wolfram Mathematica 
or its numerical eigenvalue computation routines,
which offer arbitrarily high precision (such as $10^{-100}$), 
at the expense of longer computational time.
After the energy of both true and false positives from the first phase
is computed either symbolically or with much higher precision in Wolfram Mathematica,
only true positives are left in the final set.
The programming code used to perform these two-phase searches can be found in~\cite{GitHubStuff}.

\subsection{Borderenergetic graphs of order \texorpdfstring{$12$}{12}}\label{subsc_border_12}

Furtula and Gutman \cite{FurGut2017} previously performed a computer-assisted search 
incorrectly deducing the existence of $572$ connected borderenergetic graphs of order~$12$. 
Some of the obtained graphs were not filtered properly for numerical precision errors, 
hence several of these graphs are not actually borderenergetic. 

In the first phase of our search,
we went through the set of all graphs of order~$12$, 
generated by the program~\emph{geng} from the package~\emph{nauty} \cite{McKayPip2014},
and computed their energy with precision $10^{-6}$, 
resulting in the intermediate set of $1826$ candidate graphs
consisting of true and false positives with the energy in the interval $[22-10^{-6}, 22+10^{6}]$.

In the second phase, 
we computed the energy of each candidate graph with the precision of~$10^{-100}$ using Wolfram Mathematica.
This resulted in the final list of $566$~graphs.

As an even more superior way to handling numerical errors,
we used Wolfram Mathematica to also symbolically verify that 
all the reported $566$ true positives truly are borderenergetic, thus completing the search.
As a side note, we would like to point out that 
Wolfram Mathematica does not execute expression simplification per se while performing the last equality check.
Hence, without explicitly calling \texttt{Simplify} function,
some of the borderenergetic graphs may end up missing. 
(While one could alternatively also use \texttt{FullSimplify} function, 
 this would take too much time, without any additional gain.)
Therefore, we used \texttt{Simplify} function in this extra step to confirm Theorem~\ref{borderenergetic_12_th},
and obtain the correct entry for order~12 in Table~\ref{table:borderenergetic_count}.

Among the $566$ borderenergetic graphs of order~$12$, there are seven disconnected graphs.
They all have exactly two components---six of them have an isolated vertex, 
while one graph contains $K_2$ as a component.

Finally, by comparing the resulting $559$ borderenergetic connected graphs of order~$12$ 
with the $572$ graphs found in \cite{FurGut2017}, 
it becomes apparent that there are $13$ graphs from the latter list that are not actually borderenergetic. 
These graphs, which all have energy that differs from $\mathcal{E}(K_{12})$ by less than $10^{-10}$, 
are listed in Table \ref{tab:overflow}.

\afterpage{
\begin{longtable}{ccc}
\includegraphics[width=0.160\textwidth]{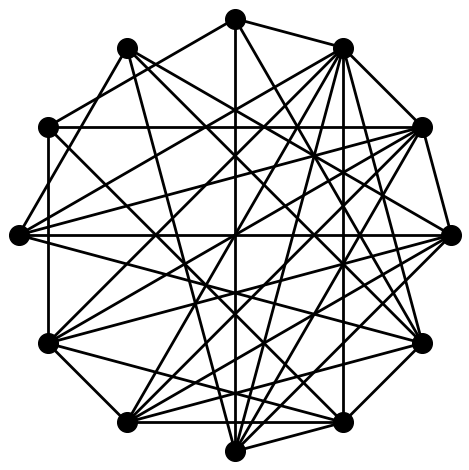} &
\includegraphics[width=0.160\textwidth]{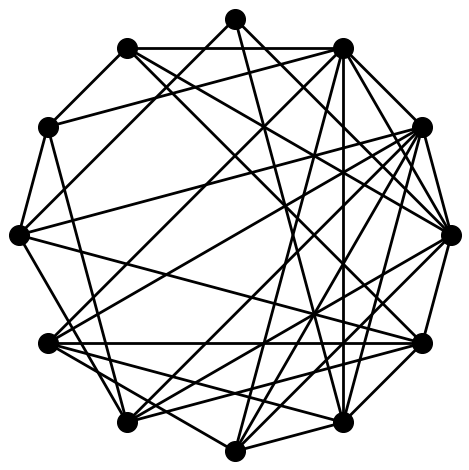} &
\includegraphics[width=0.160\textwidth]{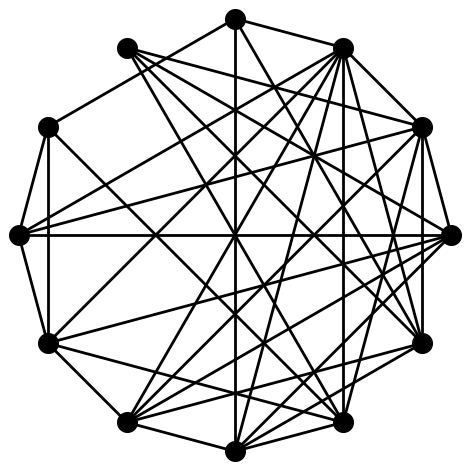} \\
\texttt{KQGM@\^{}TV`\{r|} & 
\texttt{KL?oST\textbackslash pg\}lZ} &
\texttt{KPKKJVEVajr|} \\
{\footnotesize 21.999999999908643447} &
{\footnotesize 21.999999999922612155} &
{\footnotesize 21.999999999908643447} \\[8pt]
\includegraphics[width=0.160\textwidth]{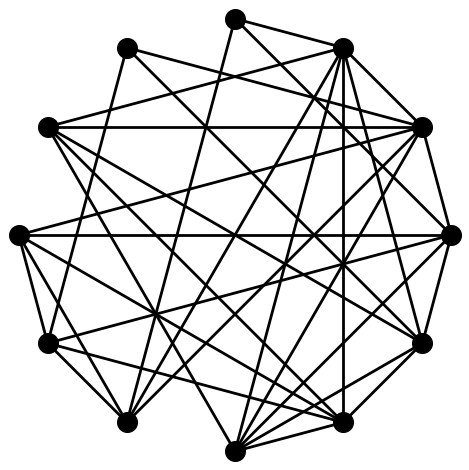} &
\includegraphics[width=0.160\textwidth]{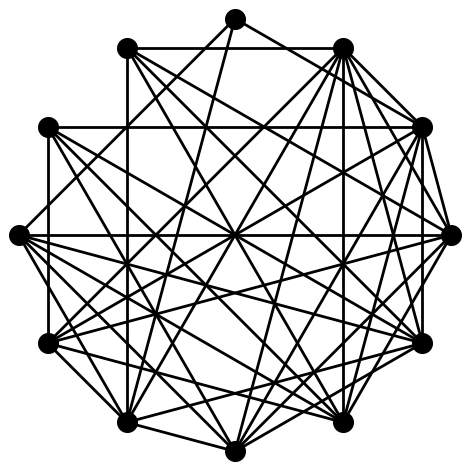} &
\includegraphics[width=0.160\textwidth]{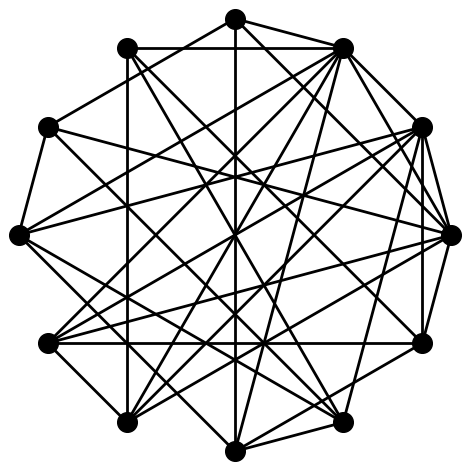} \\
\texttt{K?UX@tbejks|} &
\texttt{KCJXjpuUtVh\textasciitilde{}} & 
\texttt{KP@KbdIjGzzr} \\
{\footnotesize 21.999999999922612155} &
{\footnotesize 21.999999999930439616} &
{\footnotesize 21.999999999922612155} \\[8pt]
\includegraphics[width=0.160\textwidth]{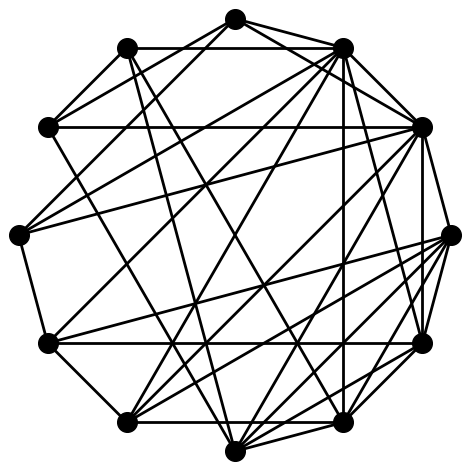} &
\includegraphics[width=0.160\textwidth]{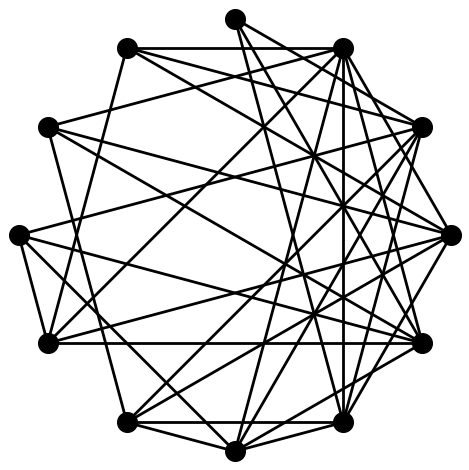} &
\includegraphics[width=0.160\textwidth]{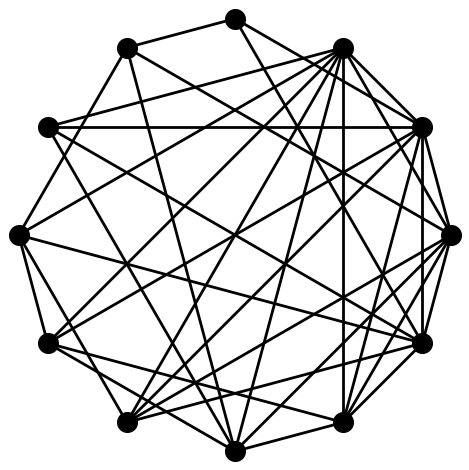} \\
\texttt{K[CJAKJB|lzl} &
\texttt{K?S\_kMyZUmL]} & 
\texttt{KaCROUtP|Zfz} \\
{\footnotesize 22.000000000040808381} &
{\footnotesize 21.999999999907887937} &
{\footnotesize 22.000000000015117356} \\[8pt]
\includegraphics[width=0.160\textwidth]{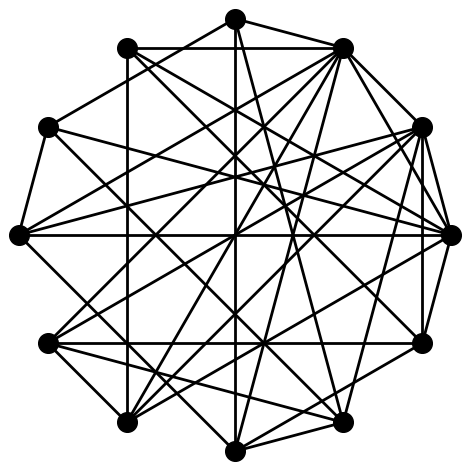} &
\includegraphics[width=0.160\textwidth]{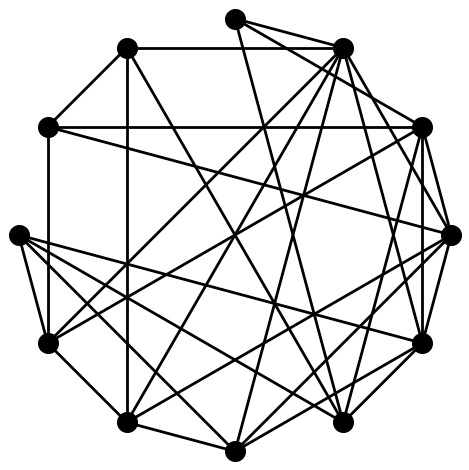} &
\includegraphics[width=0.160\textwidth]{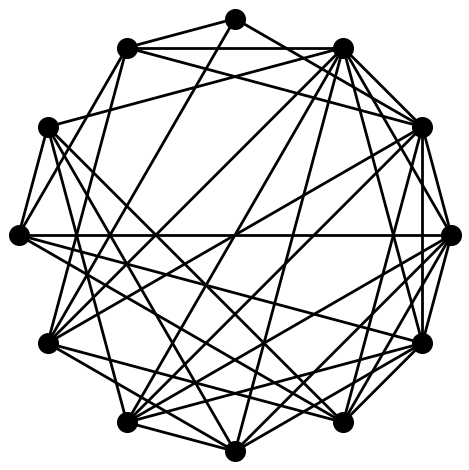} \\
\texttt{KP@KdTI\textbackslash Gzzr} &
\texttt{KGLGm\_RHlRxu} & 
\texttt{Kbo`XoVD\}Zlv} \\
{\footnotesize 21.999999999922612155} &
{\footnotesize 21.999999999967741356} &
{\footnotesize 21.999999999983254373} \\[8pt]
& \includegraphics[width=0.160\textwidth]{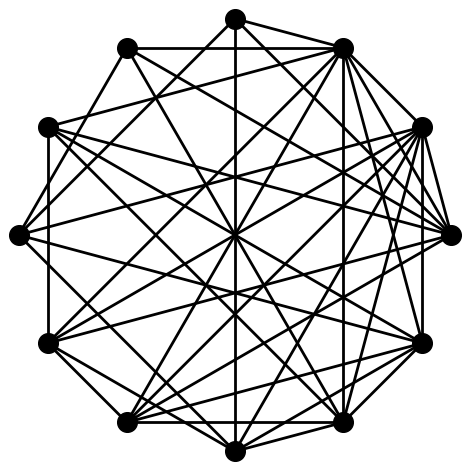} & \\
& \texttt{KEGKrKvz?\textasciitilde{}|n} & \\
& {\footnotesize 21.999999999944035209} & \\
\caption{Non-borderenergetic graphs on $12$~vertices whose energy differs from~$\mathcal{E}(K_{12})$ by less than $10^{-10}$.}
\label{tab:overflow}
\end{longtable}
}

\subsection{Borderenergetic chemical graphs on at most \texorpdfstring{$19$}{19} vertices}

%Theorems \ref{li_th_1} and \ref{li_th_2} together imply that any borderenergetic chemical graph must have a maximum degree of four. Furthermore, such a graph $G$ cannot have more than $21$ vertices, while its number of edges $|E(G)|$ is surely equal to either $2|V(G)|$ or $2|V(G)| - 1$. With this in mind, we can prove Theorem~\ref{main_chemical_th} by performing a computer-assisted search that finds all the borderenergetic chemical graphs on at most $19$ vertices, and then by showing that there exists no borderenergetic chemical graph on $20$ or $21$ vertices, as done so in Section~\ref{sc_nonexistence}.

%Taking into consideration the earlier results \cite{GongLiXuGutFur2015, LiWeiGong2015, ShaoDeng2016} together with the new search results obtained in Subsection \ref{subsc_border_12}, we may directly conclude that the only borderenergetic chemical graphs of order up to $12$ are the two already known graphs from Figure \ref{fig:chemical9}. 

Here we searched for borderenergetic graphs among all graphs of maximum degree four, 
order $n$ and size $2n - 1$ or $2n$, where $n \le 19$. 
All such graphs were again generated using \emph{geng}~\cite{McKayPip2014}.
For each $n\leq 17$ it was possible to put the whole set of generated graphs in a relatively small file 
and perform the search on a single personal computer at once.

However, 
the numbers of these graphs with $n=18$ and $n=19$ vertices 
were too large for generated files to easily fit easily on a personal computer.
Therefore,
\emph{geng} was instructed to partition the generation of graphs 
into 100 groups for $n=18$ and
into 100,000 groups for $n=19$.
These groups were then split among seven relatively modern personal computers,
which took about a month to go through the first phase of the search and
select the candidate graphs whose energy differs from~$\mathcal{E}(K_n)$ less than $10^{-8}$.
These candidates were then in the second phase checked with Wolfram Mathematica
by computing their energies with the 100-digit precision.

Besides the two borderenergetic chemical graphs on nine vertices (Fig.~\ref{fig:chemical9}),
that were already found in~\cite{LiWeiZhu2017},
this search yielded one more borderenergetic chemical graph of order~$15$, shown in Fig.~\ref{fig:chemical15}.
Therefore, these three graphs are the only three borderenergetic chemical graphs of order at most~$19$.

\section{Nonexistence of borderenergetic chemical graphs on \texorpdfstring{$20$ and $21$}{20 and 21} vertices}\label{sc_nonexistence}

In this section, we finalize the proof of Theorem \ref{main_chemical_th} by showing 
that there is no borderenergetic chemical graph on $20$ or $21$ vertices. 
Bearing in mind Theorems \ref{li_th_1} and \ref{li_th_2}, 
we may accomplish this by proving the following two theorems.

\begin{theorem}\label{berkovitz_th_1}
    For any $4$-regular graph $G$ on $n = 20$ or $21$ vertices, we have $\mathcal{E}(G) < 2(n-1)$.
\end{theorem}

\begin{theorem}\label{berkovitz_th_2}
    Let $G$ be a graph on $n=20$ or $21$ vertices and $2n - 1$ edges, such that its maximum vertex degree is at most four. Then $\mathcal{E}(G) < 2(n - 1)$.
\end{theorem}

Throughout the section, 
we will heavily rely on the next result from mathematical optimization theory, 
whose proof can be found, for example, in \cite[Theorem 5.1]{Berkovitz}.

\begin{theorem}\label{berko_th}
    Let $\bm{x}^* \in \mathbb{R}^n$ be a solution to the optimization problem
    \begin{align*}
        \min \quad & f(\bm{x}),\\
        \mathrm{s.t.} \quad & \bm{g}(\bm{x}) \le \bm{0},\\
        & \bm{h}(\bm{x}) = \bm{0},
    \end{align*}
    over an open convex set $X_0 \subseteq \mathbb{R}^n$, where $f, \bm{g}, \bm{h}$ are $C^1$ functions with domain $X_0$ and ranges in $\mathbb{R}^1, \mathbb{R}^m, \mathbb{R}^k$, respectively. Then there exists a real number $\lambda_0 \ge 0$, a vector $\bm{\lambda} \ge \bm{0}$ in $\mathbb{R}^m$, and a vector $\bm{\mu}$ in $\mathbb{R}^k$ such that
    \begin{enumerate}[label=\textbf{(\roman*)}]
        \item $(\lambda_0, \bm{\lambda}, \bm{\mu}) \neq \bm{0}$,
        \item $\langle \bm{\lambda}, \bm{g}(\bm{x}^*) \rangle = 0$, and
        \item $\lambda_0 \nabla f(\bm{x}^*) + \bm{\lambda}^\intercal \nabla \bm{g}(\bm{x}^*) + \bm{\mu}^\intercal \nabla \bm{h}(\bm{x}^*) = \bm{0}$.
    \end{enumerate}
\end{theorem}

We will also need the following short lemma, as well as two of its immediate corollaries.

\begin{lemma}\label{cool_lemma}
    The number of closed walks of length four in any graph of order~$n$ and size~$m$ is at least $\frac{8m^2}{n} - 2m$.
\end{lemma}
\begin{proof}
    Let $G$ be the given graph of order $n$ and size $m$, so that its vertex degrees are $d_1, d_2, \ldots, d_n$. 
    Let $Z_1, Z_2, Z_3$ denote the following three sets of closed walks of length four:
    \begin{align*}
        Z_1 &= \{ v_0 v_1 v_0 v_2 v_0 \colon v_0, v_1, v_2 \in V(G), v_1, v_2 \in N(v_0) \},\\
        Z_2 &= \{ v_0 v_1 v_2 v_1 v_0 \colon v_0, v_1, v_2 \in V(G), v_0, v_2 \in N(v_1) \},\\
        Z_3 &= \{ v_0 v_1 v_0 v_1 v_0 \colon v_0, v_1 \in V(G), \, v_1 \in N(v_0) \} ,
    \end{align*}
    The number of closed walks of length four in~$G$ is then at least $|Z_1| + |Z_2| - |Z_3|$. 
    Bearing in mind that
    \[
        |Z_1| = |Z_2| = \sum_{i = 1}^n d_i^2 \qquad \mbox{and} \qquad |Z_3| = \sum_{i = 1}^n d_i = 2m,
    \]
    the number of closed walks of length four is at least 
    $$
    2\sum_{i = 1}^n d_i^2 - 2m \geq 2\frac{\left(\sum_{i=1}^n d_i\right)^2}n - 2m = \frac{8m^2}n - 2m,
    $$
    by the Cauchy-Schwarz inequality.
\end{proof}

\begin{corollary}\label{cool_cor_1}
    If $G$ is a $d$-regular graph on $n$ vertices, then its number of closed walks of length four is at least $n(2d^2 - d)$.
\end{corollary}
\begin{proof}
    Since $m=\frac{nd}{2}$, 
    from Lemma~\ref{cool_lemma} the number of closed walks of length four is at least
    $\frac 8n \left(\frac{nd}{2}\right)^2 - 2 \, \frac{nd}{2} = n(2d^2 - d)$.
%    \[
%        \pushQED{\qed}
%        \frac{8 \left(\frac{nd}{2}\right)^2}{n} - 2 \, \frac{nd}{2} = \frac{2n^2 d^2}{n} - nd = n(2d^2 - d).\qedhere
%    \]
\end{proof}

\begin{corollary}\label{cool_cor_2}
    If $G$ is a graph of order $n$ and size $\frac{nd}{2} - 1$, then its number of closed walks of length four is at least $n(2d^2 - d) - 8d + 2$.
\end{corollary}
\begin{proof}
    Lemma \ref{cool_lemma} directly implies that the number of closed walks of length four is at least
    \begin{align*}
        \pushQED{\qed}
        \frac{8 \left( \frac{nd}{2} - 1\right)^2}{n} - 2 \left( \frac{nd}{2} - 1 \right) &= \frac{2 n^2 d^2 - 8nd + 8}{n} - nd + 2\\
        &> \left( 2nd^2 - 8d \right) - nd + 2\\
        &= n(2d^2 - d) - 8d + 2. \qedhere
    \end{align*}
\end{proof}

Besides that, we will use two more auxiliary lemmas. The proof of Lemma~\ref{system_lemma_1} is somewhat technical and can be found in \ref{appendix_proof}, while we choose to omit the proof of Lemma \ref{system_lemma_2} since it is entirely analogous to that of Lemma \ref{system_lemma_1}.

\begin{lemma}\label{system_lemma_1}
    For fixed parameters $\alpha, \beta \in \mathbb{N}$ and $n \in \mathbb{N}, \, n \ge 14$, such that $\alpha + \beta = n - 1$, the system of equations
    \begin{align}
        \label{oaux_1}\alpha A + \beta B &= 4n - 16,\\
        \label{oaux_2}\alpha A^2 + \beta B^2 &= 28n - 256,
    \end{align}
    in $A, B \in \mathbb{R}$ has exactly the two solutions
    \begin{equation}\label{oaux_4}
        (A, B) = \left( \frac{4n - 16}{n-1} + \frac{\sqrt{12 \beta n (n-13)}}{(n-1) \sqrt{\alpha}}, \frac{4n - 16}{n-1} - \frac{\sqrt{12 \alpha n (n-13)}}{(n-1) \sqrt{\beta}} \right)
    \end{equation}
    and
    \begin{equation}\label{oaux_5}
        (A, B) = \left( \frac{4n - 16}{n-1} - \frac{\sqrt{12 \beta n (n-13)}}{(n-1) \sqrt{\alpha}}, \frac{4n - 16}{n-1} + \frac{\sqrt{12 \alpha n (n-13)}}{(n-1) \sqrt{\beta}} \right) .
    \end{equation}
\end{lemma}

\begin{lemma}\label{system_lemma_2}
    For fixed parameters $\alpha, \beta \in \mathbb{N}$ and $n \in \mathbb{N}, \, n \ge 2$, such that $\alpha + \beta = n$, the system of equations
    \begin{align}
        \label{oaux_8}\alpha A + \beta B &= 4n - 2,\\
        \label{oaux_9}\alpha A^2 + \beta B^2 &= 28n - 30,
    \end{align}
    in $A, B \in \mathbb{R}$ has exactly the two solutions
    \begin{equation}\label{oaux_10}
        (A, B) = \left( \frac{4n - 2}{n} + \frac{\sqrt{2 \beta (6n^2 - 7n - 2)}}{n \sqrt{\alpha}}, \frac{4n - 2}{n} - \frac{\sqrt{2 \alpha (6n^2 - 7n - 2)}}{n \sqrt{\beta}} \right)
    \end{equation}
    and
    \begin{equation}\label{oaux_11}
        (A, B) = \left( \frac{4n - 2}{n} - \frac{\sqrt{2 \beta (6n^2 - 7n - 2)}}{n \sqrt{\alpha}}, \frac{4n - 2}{n} + \frac{\sqrt{2 \alpha (6n^2 - 7n - 2)}}{n \sqrt{\beta}} \right) .
    \end{equation}
\end{lemma}

We are finally in the position to provide the proofs of Theorems \ref{berkovitz_th_1} and \ref{berkovitz_th_2}.

\bigskip\noindent
\emph{Proof of Theorem \ref{berkovitz_th_1}}.\quad Let $\lambda_1 \ge \lambda_2 \ge \cdots \ge \lambda_n$ be the eigenvalues of $G$. From the Perron--Frobenius theorem, it follows that $\lambda_1 = 4$, so we immediately get
\[
    \sum_{i = 1}^n \lambda_i^2 = 4n, \quad \mbox{i.e.,} \quad \sum_{i = 2}^n \lambda_i^2 = 4n - 4^2,
\]
and by Corollary \ref{cool_cor_1}, we may also conclude that
\[
    \sum_{i = 1}^n \lambda_i^4 \ge 28n, \quad \mbox{i.e.,} \quad \sum_{i = 2}^n \lambda_i^4 \ge 28n - 4^4 .
\]
Therefore, in order to prove that $\mathcal{E}(G) < 2(n - 1)$ surely holds, it is sufficient to demonstrate that any (potential) solution to the minimization problem
\begin{align*}
    \min \quad & -\sum_{i = 1}^{n-1} x_i,\\
    \mathrm{s.t.} \quad & 28n - 256 - \sum_{i = 1}^{n-1} x_i^4 \le 0,\\
    & \sum_{i = 1}^{n-1} x_i^2 - 4n + 16 = 0,
\end{align*}
in $(x_1, x_2, \ldots, x_{n-1})$ over $\mathbb{R}^{n-1}$ satisfies $\sum_{i = 1}^{n-1} x_i < 2n - 6$. 

Note that from the above, the point $(|\lambda_2|,\dots ,|\lambda_n|)$ is feasible, 
hence the set of feasible points is nonempty, closed and bounded 
as it is contained in $[-\sqrt{4n-16},\sqrt{4n-16}]^{n-1}$. 
Since a continuous function has both a minimizer and a maximizer on a compact set, 
we get that the optimization problem has a solution, and 
we also deduce that $\sum_{i=1}^{n-1}|\lambda_i|<2n-6$.

%This is because the given optimization problem certainly has a solution. Indeed, its set of feasible points is closed, hence its boundedness by $\left[-\sqrt{4n - 16}, \sqrt{4n-16} \right]^{n-1}$ guarantees that it must be compact. Furthermore, the vector $(\lambda_2, \lambda_3, \ldots, \lambda_n)$ is clearly feasible. The desired conclusion promptly follows from the well-known fact that any real continuous function necessarily has both a minimizer and a maximizer on any nonempty compact set.

Now, let $\bm{x}^* = (x^*_1, x^*_2, \ldots, x^*_{n-1}) \in \mathbb{R}^{n-1}$ be a solution to the given minimization problem. 
By Theorem \ref{berko_th}, we have that there must exist real numbers $\xi_0, \xi_1 \ge 0$ and $\xi_2 \in \mathbb{R}$ such that:
\begin{enumerate}[label=\textbf{(\roman*)}]
    \item $(\xi_0, \xi_1, \xi_2) \neq (0, 0, 0)$;
    \item $\xi_1 = 0$ or $\sum_{i = 1}^{n-1} \left(x^*_i\right)^4 = 28n - 256$;
    \item $-\xi_0 - 4\xi_1 \left(x^*_i\right)^3 + 2 \xi_2 x^*_i = 0$ for every $1 \le i \le n-1$.
\end{enumerate}
If $\xi_1 = \xi_2 = 0$, then condition \textbf{(iii)} implies that $\xi_0 = 0$, which is impossible due to condition~\textbf{(i)}. On the other hand, if $\xi_1 = 0, \, \xi_2 \neq 0$ (resp.\ $\xi_2 = 0, \linebreak \xi_1 \neq 0$), then the strict monotonicity of $t \mapsto -\xi_0 + 2 \xi_2 t$ (resp.\ $t \mapsto -\xi_0 - 4\xi_1 t^3$) guarantees that $x^*_1 = x^*_2 = \cdots = x^*_{n-1}$. However, in this case, it follows that $x^*_i = \pm \sqrt{\frac{4n-16}{n-1}}$ for each $1 \le i \le n-1$, hence
\begin{align*}
    \sum_{i = 1}^{n-1} \left(x^*_i\right)^4 &= (n-1) \, \frac{(4n-16)^2}{(n-1)^2} = \frac{(4n-16)^2}{n-1}\\
    &= \frac{(28n - 256)(n-1) - 12n(n-13)}{n-1} < 28n - 256 .
\end{align*}
Thus, the given point is infeasible, which means that $\xi_1, \xi_2 \neq 0$ necessarily holds.

Moving on, if $\xi_0 = 0$, then from condition \textbf{(iii)} we obtain $x^*_i = 0$ or $\left(x^*_i\right)^2 = \frac{\xi_2}{2 \xi_1}$, for each $1 \le i \le n-1$. Now, if we let $\alpha \in \mathbb{N}$ denote the number of nonzero values among $x^*_1, x^*_2, \ldots, x^*_{n-1}$, then we may directly conclude that all $\alpha$ of these numbers need to have the same absolute value of $\sqrt{\frac{4n-16}{\alpha}}$. If $\alpha = n - 1$, then this directly leads to an infeasible point, as already shown, hence we may assume that $1 \le \alpha \le n - 2$. With this in mind, it follows that
\begin{align*}
    \sum_{i = 1}^{n-1} x^*_i &\le \alpha \sqrt{\frac{4n-16}{\alpha}} = \sqrt{\alpha(4n - 16)}\\
    &\le \sqrt{(n-2)(4n-16)} = \sqrt{(2n - 6)^2 - 4} < 2n - 6 ,
\end{align*}
hence any such potential minimizer would necessarily satisfy the inequality $\sum_{i = 1}^{n-1} x^*_i < 2n - 6$.

The last remaining scenario is when $\xi_0 \neq 0$, in which case we may assume without loss of generality that $\xi_0 = 1$, $\xi_1 > 0$ and $\xi_2 \neq 0$. Here, condition \textbf{(iii)} tells us that each of the $x^*_1, x^*_2, \ldots, x^*_{n-1}$ numbers represents a root of the $\mathbb{R}[t]$ polynomial $P(t) = -4\xi_1 t^3 + 2\xi_2 t - 1$. 
Vieta's formulas now yield that the product of the three roots of $P(t)$ equals $-\frac{1}{4\xi_1} < 0$, thus implying that $P(t)$ can have at most two distinct nonnegative roots. However, it is trivial to observe that all the $x^*_1, x^*_2, \ldots, x^*_{n-1}$ numbers must be nonnegative, which means that these numbers possess at most two distinct values. We have already shown that $x^*_1 = x^*_2 = \cdots = x^*_{n-1}$ leads to an infeasible point, hence we may conclude that there exist $A, B \ge 0, \, A \neq B$, such that precisely $\alpha \in \mathbb{N}, \, \alpha \le n - 2$, numbers among $x^*_1, x^*_2, \ldots, x^*_{n-1}$ are equal to $\sqrt{A}$, while the remaining $\beta = (n-1) - \alpha$ are equal to $\sqrt{B}$. Since $\xi_1 \neq 0$, condition \textbf{(ii)} allows us to reach the system of equations
\begin{equation}\label{oaux_6}
    \alpha A + \beta B = 4n - 16, \qquad \alpha A^2 + \beta B^2 = 28n - 256 .
\end{equation}

By fixing the values $\alpha, \beta \in \mathbb{N}$, it is possible to regard \eqref{oaux_6} as a system of equations in $(A, B)$ over $[0, +\infty)^2$. Due to Lemma \ref{system_lemma_1}, we know that such a system can have at most two solutions, thus making it convenient to finalize the proof of the theorem by repeating the next two steps for all the $(\alpha, \beta) = (1, n-2), (2, n-3), \ldots, \left(\lfloor \frac{n-1}{2} \rfloor, \lceil \frac{n-1}{2} \rceil \right)$ where $n\in\{20, 21\}$:
\begin{enumerate}[label=\textbf{(\arabic*)}]
    \item Apply Expressions \eqref{oaux_4} and \eqref{oaux_5} in order to find the two distinct solutions $(A_1, B_1)$ and $(A_2, B_2)$ to the system of Equations \eqref{oaux_1} and \eqref{oaux_2} over $\mathbb{R}^2$.
    \item For $i = 1, 2$, check if $(A_i, B_i)$ is indeed a solution to the given system of equations over $[0, +\infty)^2$, i.e., whether $A_i, B_i \ge 0$ is satisfied. In the case that $A_i, B_i \ge 0$ holds, verify that $\sum_{i = 1}^{n-1} x^*_i = \alpha \sqrt{A} + \beta \sqrt{B}$ is truly below $2n - 6$.
\end{enumerate}
The required verification can be done with the help of the Python script whose code is provided in \cite{GitHubStuff}. The computational results can be found in Tables \ref{tab:python_verification_1} and \ref{tab:python_verification_2} from \ref{sc_python_verification} and they show that $\sum_{i = 1}^{n-1} x^*_i < 2n - 6$ does necessarily hold for any (potential) minimizer $\bm{x}^*$. \hfill\qed

\bigskip\noindent
\emph{Proof of Theorem \ref{berkovitz_th_2}}.\quad If we let $\lambda_1 \ge \lambda_2 \ge \cdots \ge \lambda_n$ be the eigenvalues of $G$, we then have $\sum_{i = 1}^n \lambda_i^2 = 4n - 2$. Moreover, the direct implementation of Corollary \ref{cool_cor_2} yields $\sum_{i = 1}^n \lambda_i^4 \ge 28n - 30$. Thus, in order to prove that $\mathcal{E}(G) < 2(n-1)$ holds, it is enough to show that any (potential) solution to the minimization problem
\begin{align*}
    \min \quad & -\sum_{i = 1}^{n} x_i,\\
    \mathrm{s.t.} \quad & 28n - 30 - \sum_{i = 1}^{n} x_i^4 \le 0,\\
    & \sum_{i = 1}^{n} x_i^2 - 4n + 2 = 0,
\end{align*}
in $(x_1, x_2, \ldots, x_{n})$ over $\mathbb{R}^{n}$ satisfies $\sum_{i = 1}^{n} x_i < 2n - 2$. Of course, this is because the given optimization problem necessarily has a solution, which can be shown analogously as in Theorem \ref{berkovitz_th_1}.

Now, if we let $\bm{x}^* = (x^*_1, x^*_2, \ldots, x^*_n) \in \mathbb{R}^n$ be a solution to the given optimization problem, then Theorem \ref{berko_th} implies the existence of real numbers $\xi_0, \xi_1 \ge 0$ and $\xi_2 \in \mathbb{R}$ such that:
\begin{enumerate}[label=\textbf{(\roman*)}]
    \item $(\xi_0, \xi_1, \xi_2) \neq (0, 0, 0)$;
    \item $\xi_1 = 0$ or $\sum_{i = 1}^{n} \left(x^*_i\right)^4 = 28n - 30$;
    \item $-\xi_0 - 4\xi_1 \left(x^*_i\right)^3 + 2 \xi_2 x^*_i = 0$ for every $1 \le i \le n$.
\end{enumerate}
In an identical manner as in the proof of Theorem \ref{berkovitz_th_1}, it is possible to show that, unless $\xi_1, \xi_2 \neq 0$, all the $x^*_1, x^*_2, \ldots, x^*_n$ values are certainly equal. However, in this case, we promptly obtain $x^*_i = \pm \sqrt{\frac{4n-2}{n}}$ for each $1 \le i \le n$, hence
\begin{align*}
    \sum_{i = 1}^{n} \left(x^*_i\right)^4 &= n \, \frac{(4n-2)^2}{n^2} = \frac{(4n-2)^2}{n}\\
    &= \frac{n(28n - 30) - (12n^2 - 14n - 4)}{n} < 28n - 30 .
\end{align*}
Therefore, such a point would be infeasible, which is impossible. This guarantees that $\xi_1, \xi_2 \neq 0$ must be true.

We will now demonstrate that the numbers $x^*_1, x^*_2, \ldots, x^*_n$ necessarily possess exactly two distinct values. 
To begin with, it is obvious that all of these numbers are surely nonnegative. Furthermore, if $\xi_0 = 0$, then condition \textbf{(iii)} dictates that either $x^*_i = 0$ or $\left(x^*_i\right)^2 = \frac{\xi_2}{2\xi_1}$, for each $1 \le i \le n$. Since it is impossible for all the $x^*_1, x^*_2, \ldots, x^*_n$ numbers to be equal, as previously shown, the nonnegativity of these numbers promptly implies that they must have precisely two distinct values. On the other hand, if $\xi_0 \neq 0$, then we may assume without loss of generality that $\xi_0 = 1$, $\xi_1 > 0$ and $\xi_2 \neq 0$. However, in this case, it is possible to use the same approach as in the proof of Theorem \ref{berkovitz_th_1} in order to observe that the $x^*_1, x^*_2, \ldots, x^*_n$ numbers surely possess two distinct values once again. Taking everything into consideration, we conclude that there exist $A, B \ge 0, \, A \neq B$, such that exactly $\alpha \in \mathbb{N}, \, \alpha \le n - 1$, numbers among $x^*_1, x^*_2, \ldots, x^*_n$ are equal to $\sqrt{A}$, while the remaining $\beta = n - \alpha$ are equal to $\sqrt{B}$. Since $\xi_1 \neq 0$, condition \textbf{(ii)} implies
\begin{equation}\label{oaux_7}
    \alpha A + \beta B = 4n - 2, \qquad \alpha A^2 + \beta B^2 = 28n - 30 .
\end{equation}

By fixing the values $\alpha, \beta \in \mathbb{N}$, it is convenient to regard \eqref{oaux_7} as a system of equations in $(A, B)$ over $[0, +\infty)^2$. From Lemma \ref{system_lemma_2}, it follows that such a system can have at most two solutions, hence we can complete the proof of the theorem by repeating the following two steps for all the $(\alpha, \beta) = (1, n-1), (2, n-2), \ldots, \left(\lfloor \frac{n}{2} \rfloor, \lceil \frac{n}{2} \rceil \right)$ where $n\in\{20, 21\}$:
\begin{enumerate}[label=\textbf{(\arabic*)}]
    \item Use Expressions \eqref{oaux_10} and \eqref{oaux_11} to find the two distinct solutions $(A_1, B_1)$ and $(A_2, B_2)$ to the system of Equations \eqref{oaux_8} and \eqref{oaux_9} over $\mathbb{R}^2$.
    \item For $i = 1, 2$, check if $(A_i, B_i)$ is truly a solution to the given system of equations over $[0, +\infty)^2$, i.e., whether $A_i, B_i \ge 0$ holds. In the case that $A_i, B_i \ge 0$ is satisfied, verify that $\sum_{i = 1}^{n} x^*_i = \alpha \sqrt{A} + \beta \sqrt{B}$ is indeed below $2n - 2$.
\end{enumerate}
As in the proof of Theorem \ref{berkovitz_th_1}, the required computation and verification can be performed via the Python script whose code is given in \cite{GitHubStuff}. From the computational results given in Tables \ref{tab:python_verification_3} and \ref{tab:python_verification_4} from \ref{sc_python_verification} we may immediately conclude that $\sum_{i=1}^n x^*_i < 2n - 2$ does hold for any (potential) minimizer $\bm{x}^*$. \hfill\qed

\section{Conclusions}\label{sc_conclusion}

The main results of this paper are the full classifications of 
the borderenergetic chemical graphs and the borderenergetic graphs of order~$12$.
A computer-assisted search was used to find all the borderenergetic graphs on $12$ vertices, 
alongside the borderenergetic chemical graphs of order at most $19$. 
The basic idea behind the search was to divide it into two phases: coarse sieve and fine sieve.
In the coarse sieve phase,
graph energy of the whole set of graphs was computed with standard numerical precision 
(in some cases computations were split among multiple computers running in parallel on smaller subsets),
and lower precision threshold was used to select relatively few candidate graphs 
whose reported energy does not differ ``too much'' from $\mathcal{E}(K_n)$,
and which include both true and false positives.
In the fine sieve phase, 
energy of all candidate graphs was recomputed with a 100-digit precision 
using arbitrary precision routines in Wolfram Mathematica
in order to separate true positive from false positives,
who simply happen to have the same energy as $\mathcal{E}(K_n)$
on the first ten or so digits.
Wolfram Mathematica also enables one to compute the energy symbolically,
which can be used as an alternative second phase to high precision numerical computations.

The given search technique makes it possible to find all the borderenergetic graphs of order~$12$, 
which corrects the previously obtained results by Furtula and Gutman \cite{FurGut2017}. 
It also enables one to find all three borderenergetic chemical graphs of order up to~$19$.
While chemical graphs on 20 and 21 vertices were too numerous for an exhaustive search,
it was possible to disprove the existence of borderenergetic chemical graphs on 20 and 21 vertices 
using a different technique based on optimization theory. 
Consequently we obtained Theorem~\ref{main_chemical_th}, 
which completes the search for borderenergetic chemical graphs initiated by Li, Wei and Zhu \cite{LiWeiZhu2017}.

\section*{Acknowledgements}

PC is supported by the MTA-R\'enyi Counting in Sparse Graphs ``Momentum'' Research Group. 
PC is also supported by the Dynasnet European Research Council Synergy project---grant number ERC-2018-SYG 810115.
ID is supported by the Science Fund of the Republic of Serbia, grant \#6767, 
Lazy walk counts and spectral radius of threshold graphs --- LZWK.

\section*{Conflict of interest}

The authors declare that they have no conflict of interest.

%% The Appendices part is started with the command \appendix;
%% appendix sections are then done as normal sections
\appendix
\section{Proof of Lemma \ref{system_lemma_1}}
\label{appendix_proof}

Here, we provide the technical proof of Lemma \ref{system_lemma_1}.

\bigskip\noindent
\emph{Proof of Lemma \ref{system_lemma_1}}.\quad From Equation~\eqref{oaux_1} we immediately get
\begin{equation}
    \label{oaux_3} B = \frac{(4n-16) - \alpha A}{\beta} .
\end{equation}
Now, by plugging in Equation~\eqref{oaux_3} into Equation~\eqref{oaux_2}, we further obtain the following equalities:
\begin{align*}
    \alpha A^2 + \beta \, \frac{((4n - 16) - \alpha A)^2}{\beta^2} &= 28n - 256,\\
    \alpha \beta A^2 + ((4n - 16) - \alpha A)^2 &= \beta (28n - 256),\\
    (\alpha \beta + \alpha^2) A^2 - 2 \alpha(4n - 16) A &= \beta(28n - 256) - (4n-16)^2,\\
    \alpha(n-1) A^2 - 2 \alpha(4n - 16) A &= \beta(28n - 256) - (4n-16)^2,\\
    \alpha^2 A^2 - \frac{2\alpha^2 (4n - 16)}{n - 1} A &= \frac{\alpha\beta(28n - 256) - \alpha(4n-16)^2}{n-1},
\end{align*}
which then lead us to
\begin{align*}
    &\left( \alpha A - \frac{\alpha(4n - 16)}{n-1} \right)^2 = \frac{\alpha\beta(28n - 256) - \alpha(4n-16)^2}{n-1} + \frac{\alpha^2 (4n - 16)^2}{(n-1)^2},\\
    &\hspace{2.247cm} = \frac{\alpha\beta(n-1)(28n - 256) + (4n-16)^2 (\alpha^2 - (n-1)\alpha)}{(n-1)^2},\\
    &\hspace{2.247cm} = \frac{\alpha\beta(n-1)(28n - 256) - \alpha\beta (4n-16)^2}{(n-1)^2},\\
    &\hspace{2.247cm} = \frac{\alpha\beta ((n-1)(28n - 256) - (4n-16)^2)}{(n-1)^2},\\
    &\hspace{2.247cm} = \frac{12\alpha\beta n(n-13)}{(n-1)^2} .
\end{align*}
Thus, we have that either
\[
    \alpha A = \frac{\alpha(4n - 16) + \sqrt{12 \alpha \beta n (n-13)}}{n-1}
\]
and
\[
    \beta B = \frac{\beta(4n - 16) - \sqrt{12 \alpha \beta n (n-13)}}{n-1} ,
\]
or
\[
    \alpha A = \frac{\alpha(4n - 16) - \sqrt{12 \alpha \beta n (n-13)}}{n-1}
\]
and
\[
    \beta B = \frac{\beta(4n - 16) + \sqrt{12 \alpha \beta n (n-13)}}{n-1} .
\]
The two distinct solutions follow immediately from here. \hfill\qed

\section{Python verification results for Theorems \ref{berkovitz_th_1} and \ref{berkovitz_th_2}}\label{sc_python_verification}

{\footnotesize
\begin{longtable}{cccccccc}
$\alpha$ & $\beta$ & $A_1$ & $B_1$ & $\alpha \sqrt{A_1} + \beta \sqrt{B_1}$ & $A_2$ & $B_2$ & $\alpha \sqrt{A_2} + \beta \sqrt{B_2}$\\
\hline
1 & 18 & $12.5209$ & $2.8600$ & $33.9790$ & $-5.7840$ & $3.8769$ & --- \\
2 & 17 & $9.6578$ & $2.6285$ & $33.7769$ & $-2.9210$ & $4.1084$ & --- \\
3 & 16 & $8.3504$ & $2.4343$ & $33.6327$ & $-1.6135$ & $4.3025$ & --- \\
4 & 15 & $7.5459$ & $2.2544$ & $33.5100$ & $-0.8091$ & $4.4824$ & --- \\
5 & 14 & $6.9782$ & $2.0792$ & $33.3954$ & $-0.2414$ & $4.6576$ & --- \\
6 & 13 & $6.5438$ & $1.9029$ & $33.2813$ & $0.1930$ & $4.8340$ & $31.2183$ \\
7 & 12 & $6.1929$ & $1.7208$ & $33.1614$ & $0.5439$ & $5.0161$ & $32.0384$ \\
8 & 11 & $5.8980$ & $1.5287$ & $33.0292$ & $0.8388$ & $5.2081$ & $32.4304$ \\
9 & 10 & $5.6424$ & $1.3219$ & $32.8756$ & $1.0945$ & $5.4150$ & $32.6856$ \\
\caption{Python verification results for Theorem \ref{berkovitz_th_1} and the case $n = 20$.}
\label{tab:python_verification_1}
\end{longtable}}

{\footnotesize
\begin{longtable}{cccccccc}
$\alpha$ & $\beta$ & $A_1$ & $B_1$ & $\alpha \sqrt{A_1} + \beta \sqrt{B_1}$ & $A_2$ & $B_2$ & $\alpha \sqrt{A_2} + \beta \sqrt{B_2}$\\
\hline
1 & 19 & $13.1857$ & $2.8850$ & $35.9031$ & $-6.3857$ & $3.9150$ & --- \\
2 & 18 & $10.1350$ & $2.6517$ & $35.6782$ & $-3.3350$ & $4.1483$ & --- \\
3 & 17 & $8.7442$ & $2.4569$ & $35.5179$ & $-1.9442$ & $4.3431$ & --- \\
4 & 16 & $7.8900$ & $2.2775$ & $35.3819$ & $-1.0900$ & $4.5225$ & --- \\
5 & 15 & $7.2884$ & $2.1039$ & $35.2556$ & $-0.4884$ & $4.6961$ & --- \\
6 & 14 & $6.8293$ & $1.9303$ & $35.1307$ & $-0.0293$ & $4.8697$ & --- \\
7 & 13 & $6.4594$ & $1.7526$ & $35.0010$ & $0.3406$ & $5.0474$ & $33.2915$ \\
8 & 12 & $6.1495$ & $1.5670$ & $34.8601$ & $0.6505$ & $5.2330$ & $33.9030$ \\
9 & 11 & $5.8819$ & $1.3693$ & $34.6994$ & $0.9181$ & $5.4307$ & $34.2576$ \\
10 & 10 & $5.6450$ & $1.1550$ & $34.5063$ & $1.1550$ & $5.6450$ & $34.5063$ \\
\caption{Python verification results for Theorem \ref{berkovitz_th_1} and the case $n = 21$.}
\label{tab:python_verification_2}
\end{longtable}}

{\footnotesize
\begin{longtable}{cccccccc}
$\alpha$ & $\beta$ & $A_1$ & $B_1$ & $\alpha \sqrt{A_1} + \beta \sqrt{B_1}$ & $A_2$ & $B_2$ & $\alpha \sqrt{A_2} + \beta \sqrt{B_2}$\\
\hline
1 & 19 & $18.5462$ & $3.1291$ & $37.9164$ & $-10.7462$ & $4.6709$ & --- \\
2 & 18 & $13.9802$ & $2.7800$ & $37.4899$ & $-6.1802$ & $5.0200$ & --- \\
3 & 17 & $11.8985$ & $2.4885$ & $37.1657$ & $-4.0985$ & $5.3115$ & --- \\
4 & 16 & $10.6201$ & $2.2200$ & $36.8747$ & $-2.8201$ & $5.5800$ & --- \\
5 & 15 & $9.7198$ & $1.9601$ & $36.5887$ & $-1.9198$ & $5.8399$ & --- \\
6 & 14 & $9.0326$ & $1.7003$ & $36.2881$ & $-1.2326$ & $6.0997$ & --- \\
7 & 13 & $8.4790$ & $1.4344$ & $35.9527$ & $-0.6790$ & $6.3656$ & --- \\
8 & 12 & $8.0152$ & $1.1565$ & $35.5539$ & $-0.2152$ & $6.6435$ & --- \\
9 & 11 & $7.6147$ & $0.8607$ & $35.0404$ & $0.1853$ & $6.9393$ & $32.8511$ \\
10 & 10 & $7.2601$ & $0.5399$ & $34.2926$ & $0.5399$ & $7.2601$ & $34.2926$ \\
\caption{Python verification results for Theorem \ref{berkovitz_th_2} and the case $n = 20$.}
\label{tab:python_verification_3}
\end{longtable}}

{\footnotesize
\begin{longtable}{cccccccc}
$\alpha$ & $\beta$ & $A_1$ & $B_1$ & $\alpha \sqrt{A_1} + \beta \sqrt{B_1}$ & $A_2$ & $B_2$ & $\alpha \sqrt{A_2} + \beta \sqrt{B_2}$\\
\hline
1 & 20 & $18.9542$ & $3.1523$ & $39.8630$ & $-11.1447$ & $4.6572$ & --- \\
2 & 19 & $14.2769$ & $2.8130$ & $39.4235$ & $-6.4673$ & $4.9966$ & --- \\
3 & 18 & $12.1477$ & $2.5309$ & $39.0921$ & $-4.3381$ & $5.2786$ & --- \\
4 & 17 & $10.8422$ & $2.2724$ & $38.7977$ & $-3.0327$ & $5.5371$ & --- \\
5 & 16 & $9.9245$ & $2.0236$ & $38.5120$ & $-2.1150$ & $5.7859$ & --- \\
6 & 15 & $9.2255$ & $1.7765$ & $38.2167$ & $-1.4160$ & $6.0331$ & --- \\
7 & 14 & $8.6638$ & $1.5252$ & $37.8941$ & $-0.8543$ & $6.2843$ & --- \\
8 & 13 & $8.1945$ & $1.2649$ & $37.5218$ & $-0.3850$ & $6.5446$ & --- \\
9 & 12 & $7.7905$ & $0.9905$ & $37.0629$ & $0.0190$ & $6.8191$ & $32.5771$ \\
10 & 11 & $7.4342$ & $0.6962$ & $36.4440$ & $0.3754$ & $7.1133$ & $35.4645$ \\
\caption{Python verification results for Theorem \ref{berkovitz_th_2} and the case $n = 21$.}
\label{tab:python_verification_4}
\end{longtable}}

\end{document}